\documentclass[12pt, twoside]{article} 
\usepackage{amssymb,amsfonts,amsmath,amsbsy,theorem,amscd, bbm, xspace, ifthen}

\usepackage[T1]{fontenc}
\newcommand{\dd}{\text{\it \dj\hspace{1pt}}}
\numberwithin{equation}{section}
\newcommand{\C}{\ensuremath{\mathbb{C}}}
\newcommand{\R}{\ensuremath{\mathbb{R}}}

\newcommand{\Rn}{\ensuremath{\mathbb{R}^n}}

\newcommand{\N}{\ensuremath{\mathbb{N}}}

\newcommand{\F}{\ensuremath{\mathcal{F}}}

\newcommand{\SD}{\ensuremath{\mathcal{S}}}

\newcommand{\im}{\operatorname{Im}}
\newcommand{\re}{\operatorname{Re}}

\newcommand{\sd}{\, d}
\newcommand{\supp}{\operatorname{supp}}
\newcommand{\eps}{\ensuremath{\varepsilon}}
\newcommand{\weight}[1]{\langle #1\rangle}

\newcommand{\lh}{\mathcal{C}_0}
\newcommand{\sC}{\mathcal{C}}

\newtheorem{thm}{Theorem}[section]
\newtheorem{cor}[thm]{Corollary}
\newtheorem{lem}[thm]{Lemma}
\newtheorem{defn}[thm]{Definition}
\newtheorem{theorem}[thm]{Theorem}
\newtheorem{prop}[thm]{Proposition}

\newtheorem{claim*}{Claim}

\theorembodyfont{\rmfamily}
\newtheorem{rem}[thm]{Remark}

\newenvironment{proof*}[1]{{\bf Proof
#1:}}{\hspace*{\fill}\rule{1.2ex}{1.2ex}\\ } 
\newenvironment{proof}{{\bf
Proof:\,}}{\hspace*{\fill}\rule{1.2ex}{1.2ex}\\ }

\newcommand{\supl}{\sup\limits}
\newcommand{\liml}{\lim\limits}
\newcommand{\suml}{\sum\limits}
\newcommand{\maxl}{\max\limits}

\newcommand{\infl}{\inf\limits}
\newcommand{\sB}{\mathcal{B}}
\newcommand{\sF}{\mathcal{F}}

\newcommand{\sA}{\mathcal{A}}
\newcommand{\sD}{\mathcal{D}}
\newcommand{\PP}{\mathbbm{P}}
\newcommand{\E}{\mathbbm{E}}
\newcommand{\Q}{\mathbbm{Q}}
\newcommand{\il}{\int\limits}
\setboolean{@mparswitch}{false}

\begin{document}
\reversemarginpar
\begin{titlepage}
\title{The Cauchy problem and the martingale problem for integro-differential operators with non-smooth kernels}
\author{H. Abels\footnote{Partially supported by DFG (German Science Foundation) through SFB 611} \; and M. Kassmann\footnote{Partially supported by DFG (German Science Foundation) through SFB 611}}
\end{titlepage}

\maketitle

\begin{abstract}
We consider the linear integro-differential operator $L$ defined by
\[ Lu(x) =\int_{\Rn}
\left(u(x+y)-u(x)-\mathbbm{1}_{[1,2]}(\alpha)\mathbbm{1}_{\{|y|\leq 2\}}(y)y\cdot \nabla u(x)
  \right)k(x,y) \sd y \,. \]
Here the kernel $k(x,y)$ behaves like $|y|^{-d-\alpha}$, $\alpha \in (0,2)$, for small $y$ and is H\"older-continuous in the first variable, precise definitions are given below. The aim of this work is twofold. On one hand, we study the unique solvability of the Cauchy problem corresponding to $L$. On the other hand, we study the martingale problem for $L$. The analytic results obtained for the deterministic parabolic equation guarantee that the martingale problem is well-posed. Our strategy follows the classical path of Stroock-Varadhan. The assumptions allow for cases that have not been dealt with so far. 
\end{abstract}


\noindent{\bf AMS-Classification:} 47G20, 47G30, 60J75, 60J35, 60G07, 35K99,  35B65, 47A60

\newpage

\section{Introduction}

A linear operator $A \colon C^2_0(\R^n) \to C(\R^n)$ is said to satisfy the global maximum principle if $A u(x^*) \leq 0$ for all $x^* \in \{x\in\R^n; u(x) \geq u(y) \; \forall \, y \in \R^n\}$. It is well-known that infinitesimal generators of strongly continuous contraction semi-groups on $C_0(\R^n)$ generating Markov processes satisfy the global maximum principle. Surprisingly, the global maximum principle implies already a certain structure of $A$, see \cite{Cou65}. More precisely, $A$ is the sum of a possibly degenerate elliptic diffusion operator with bounded coefficients, a drift and a jump part which we call $L$. Since $L$ alone generates pure jump processes which generalize L\'{e}vy processes it is sometimes called a L\'{e}vy-type operator, see \cite{JaSc01}, \cite{Bas04}, \cite{Jac05} and \cite{Kas06} for surveys.

It is the aim of this work to study important properties of the operator $L$ which is defined by 
\begin{equation}\label{eq:DefnL}
  Lu(x) =\int_{\Rn} \left(u(x+y)-u(x)-\mathbbm{1}_{B_2}(y)y\cdot\nabla u(x)
  \right)k(x,y) \sd y
\end{equation}
if  $1\leq \alpha <2$ and
\begin{equation}
\label{eq:DefnL2}
  Lu(x) =\int_{\Rn} \left(u(x+y)-u(x)\right)k(x,y) \sd y 
\end{equation}
if $0<\alpha<1$.
Here $k\colon \Rn\times \left(\Rn\setminus \{0\}\right)\to (0,\infty)$ is
H\"older continuous of order $\tau\in (0,1)$ in $x\in\R^n$, measurable in
$y\in\R^n\setminus\{0\}$ and can be decomposed as $k=k_1+k_2$ such that $k_1(x,y)=0$ for $|y|\geq 2$, $k_1$ is $(n+1)$-times differentiable in $y$, and the following estimates are satisfied:
\begin{alignat}{2}\label{eq:kernelEstim1}
   \|\partial_y^\beta k_1(.,y)\|_{C^\tau(\Rn)}&\leq C |y|^{-n-\alpha-|\beta|},& \qquad& 0<|y|\leq 2, \\\label{eq:kernelEstim2}
  k_1(x,y) &\geq c|y|^{-n-\alpha}, &\qquad& 0<|y|\leq 1,x\in\Rn, \\\label{eq:kernelEstim3}
  \|k_2(.,y)\|_{C^\tau(\Rn)} &\leq C|y|^{-n-\alpha'},& \qquad& 0<|y|\leq 1, \\\label{eq:kernelEstim4}
 \int_{|y|\geq 1} \|k_2(.,y)\|_{C^\tau(\Rn)} \sd y &<\infty,\\\label{eq:kernelEstim5} 
 \lim_{|y|\to \infty}\|k_2(.,y)\|_{C^\tau(\Rn)} =0
\end{alignat} 
for all $\beta\in\N_0^n$ with $|\beta|\leq N:=n+1$,
where $0\leq\alpha'< \alpha<2$. 
There are many examples satisfying these assumptions, see the discussion below. A model case is given by $k(x,y)=c |y|^{-n-\alpha}$,$y \ne 0$, which leads to $L=-(-\Delta)^{\alpha/2}$. 

Our main result concerning the Cauchy-Problem for $L$ is given by the
following theorem. In the following $\sC^s(\Rn)$, $s>0$, denotes the
H\"older-Zygmund space and $\sC^s_0(\Rn)= \overline{C_0^\infty(\Rn)}^{\|.\|_{\sC^s}}$. For a precise definition of the function spaces we refer to
Section \ref{subsec:prelims} below.
\begin{thm}\label{thm:main1}
  Let $k$ satisfy (\ref{eq:kernelEstim1})-(\ref{eq:kernelEstim3}),
  let $L$ be defined as in (\ref{eq:DefnL}), and let $T>0$, $0<s<\tau$, $0<\theta<1$. Then for every $f\in C^\theta([0,T];\lh^s(\Rn))$ with $f(0)=0$ there is a unique $u\in C^{1,\theta}([0,T];\lh^{s}(\Rn))\cap C^{\theta}([0,T];\lh^{s+\alpha}(\Rn))$ solving
  \begin{alignat}{2}\label{eq:Parabol1}
    \partial_t u - Lu &= f &\qquad& \text{in}\ (0,T)\times \Rn,\\\label{eq:Parabol2}
    u(0,\cdot)&= 0 && \text{in}\ \Rn. 
  \end{alignat}
If $f$ is non-negative, then $u$ is non-negative as well. 
\end{thm}
The latter theorem will be a direct consequence of the fact that $L$ generates an analytic semi-group on $\lh^s(\Rn)$ with $0<s<\tau$. In order to prove this we will construct an approximate resolvent to $L$ using pseudodifferential operators with non-smooth symbols.

Let us state the martingale problem.  By $\sD([0,\infty);\R^n)$ we denote the
space of all c\`{a}dl\`{a}g paths. We refer the reader to Section \ref{sec:MP}
below for a precise definition and a short discussion of
$\sD([0,\infty);\R^n)$.  A probability measure $\PP^\mu$ on
$\sD([0,\infty);\R^n)$ is said to be a solution to the martingale problem for
$(L, D(L))$ with domain $D(L)$ being contained in the set of bounded functions
$f\colon \R^n \to \R$, $L$ defined as in (\ref {eq:DefnL}) and $\mu$ a probability measure on $\R^n$ if, for any $\phi \in D(L)$   
\[ \Big( \phi(\Pi_t)-\phi(\Pi_0)- \il_0^t (L\phi)(\Pi_s) \, ds \Big)_{t \geq 0} \]
is a $\PP^\mu$-martingale with respect to the filtration $\big(\sigma(\Pi_s;
s\leq t)\big)_{t\geq 0}$ and $\PP^\mu(\Pi_0 =\mu)=1$. Here $\Pi$ is the usual
coordinate process, i.e., $\Pi\colon[0,\infty) \times \sD([0,\infty);\R^n) \to
\R^n$, $\Pi_t(\omega)=\omega(t)$. If for every $\mu$ there is a unique
solution $\PP^{\mu}$ of the martingale problem, we say that the martingale problem for $(L,D(L))$ is well-posed. 

Our main result concerning the martingale problem reads as follows. 

\begin{thm}\label{thm:main2}
Let $L$ be defined as above. Then the martingale problem for $(L,C^\infty_0(\R^n))$ is well-posed. 
\end{thm}

Studying the existence of pure jump processes, i.e., processes without a diffusion component, together with their properties is a field of still increasing interest. We list some references dealing with the martingale problem for non-local operators such as $L$. In the case $k(x,y)=k(y)$ with $k$ as in (\ref{eq:DefnL}) $L$ is a generator of a L\'{e}vy jump process, i.e., a jump process with independent stationary increments.  There are different and more elegant approaches than the martingale problem to the existence of a corresponding process, see \cite{Ber96}, \cite{Sat99}.

The martingale problem for an operator of the form $A+L$ where $A$ is a non-degenerate elliptic operator and $L$ is an operator of our type has been studied first in \cite{Kom73}, \cite{Str75}, \cite{LeMa76}. Since $A$ is a second order operator $L$ is a lower order perturbation of $A$ for many questions. \cite{Kom84a}, \cite{Kom84b} seem to be the first articles treating the martingale problem for pure jump processes generated by operators like $L$. The main assumptions are that $k(x,y)$ is a perturbation of $\widetilde{k}(x,y)=|y|^{-d-\alpha}$, $y\ne 0$, together with quite strong regularity assumptions. More general results have been obtained in \cite{NeTs89} using techniques from partial differential equations. In the latter article $k(x,y)$ is assumed to be twice continuously differentiable in the first variable. 

Strong results on the well-posedness have been obtained in \cite{MiPr92a}, \cite{MiPr92b}, \cite{MiPr93}. The authors use a setup similar to the one of the so called Calderon-Zygmund approach in the theory of partial differential equations. In \cite{MiPr92a}, \cite{MiPr92b} $k(x,y)$ is assumed to be only continuous in the first variable but some additional homogeneity is assumed in the second variable. To add a personal comment, these results have been underestimated in the literature from our point of view.  This is maybe due to the fact that the journal is not available easily and that the articles are written in a somewhat dense style.

Using pseudodifferential operators and anisotropic Sobolev spaces built with
continuous negative definite functions \cite{Hoh94} proves well-posedness of
the martingale problem under assumptions like $x\mapsto k(x,y)$ $\in
C^{3n}(\R^n)$ but allowing for a more general dependence of $k(x,y)$ on $y$. Moreover, the extension of $L$ to a generator of a Feller semi-group is discussed. See \cite{BLR99} for similar techniques in infinite dimensions and \cite{OkvC96} for related questions. In the setting of \cite{Hoh94} a parametrix for the pseudodifferential operator is constructed in \cite{Boe05}. These results do not apply to our setting since we assume only H\"{o}lder regularity of the mapping $x\mapsto k(x,y)$.

The results of \cite{NeTs89}, \cite{MiPr92a}, \cite{MiPr92b}, \cite{Hoh94} and
the ones in the present work do not imply one another but have a large region
of intersection. The assumptions on the $x$- dependence of $k(x,y)$ in
\cite{NeTs89}, \cite{Hoh94}, \cite{MiPr93} are more restrictive but the
assumptions on the $y$-dependence are partly weaker than ours. The situation is reversed when comparing our results to \cite{MiPr92a}, \cite{MiPr92b}. Our techniques solving the Cauchy problem are different from \cite{NeTs89}, \cite{Hoh94} and \cite{MiPr92a}. 

The authors of \cite{EIK04} prove solvability of the Cauchy problem
for a time dependent pseudodifferential operator $L(t)=p(t,x,D_x)$ where
the principal part of the symbol $p(t,x,\xi)$ is homogeneous in $\xi$ of degree
$\alpha\in [1,2]$ and uniformly H\"older continuous in $(t,x)$. Their results do not apply to the uniqueness for solutions of the martingale problem since sufficient regularity of solutions to the Cauchy problem is not provided. 

In the above list we do not mention results concerning what is sometimes called ''stable-like'' cases, i.e. when $k(x,y) \approx |y|^{-d-\alpha(x)}$, $y \ne 0$.  Well-posedness of the martingale problem is proved in one spatial dimension in \cite{Bas88} when $\alpha(\cdot)$ is Dini-continuous. Uniqueness problems for stochastic differential equations in similar situations but including higher dimensions and also diffusion coefficients are considered in \cite{Tsu92}. The techniques of \cite{Bas88} can be extended to higher dimensions and to a larger class of problems, see the forthcoming PhD-thesis \cite{Tang}. See \cite{JaLe93}, \cite{KiNe97} for results on the question when the linear operators of type $L$ extend to generators of Feller processes in the case when the $y$-singularity of $k(x,y)$ is of variable order. \cite{Hoh00} provides such a result together with well-posedness of the martingale problem when $x \mapsto \alpha(x)$ is smooth where $\alpha(x)$ is the order of differentiability of $L$.

One scope of this contribution is to present an application of the theory of pseudodifferential operators with non-smooth coefficients to jump processes. We hope to draw the attention of probabilists to this method. 

\section{The Cauchy problem for L\'evy-type Operators}

\subsection{Preliminaries and Notation}\label{subsec:prelims}

The characteristic function of a set $A$ is denoted by $\mathbbm{1}_A$. Furthermore,
we define $\weight{\xi}:= (1+|\xi|^2)^\frac12$ for  $\xi\in\Rn$.
Moreover, we define $\Sigma_\delta:=\{z\in\C \setminus \{0\}: |\arg z|< \delta\}$ for $0<\delta\leq\pi$.

As usual, $C_0^\infty(\Rn)$ denotes the set of all smooth and compactly
supported functions $f\colon \Rn\to \R$, $\SD(\Rn)$ denotes the space of all
smooth and rapidely decreasing functions, and $\SD'(\Rn)= (\SD(\Rn))'$ the
space of tempered distributions. $C^k(\R^n), k \in \N$, shall be the usual Banach space of continuous functions with bounded continuous derivatives up to order $k$. By $C_0^k(\R^n)$ we denote the closure of $C^\infty_0(\R^n)$ with  respect to the norm of $C^k(\R^n)$. $C^{s}(M;X)$, where $s\in (0,1)$, $M\subseteq
\Rn$, M closed, and $X$ is a Banach space,  is the space of uniformly bounded
H\"older continuous functions $f\colon M\to X$ of order $s$ with uniformly
bounded H\"older constant. Moreover, $C^s(M)=C^s(M;\R)$ and $f\in
C^{1,s}([0,T];X)$ iff $f\colon [0,T]\to X$ is continuously differentiable and
$\frac{d}{dt} f\in C^s([0,T];X)$. 
Finally, if $f\colon \Rn\to\R$, we define $(\tau_h f)(x)= f(x+h)$, $x,h\in\Rn$, and 
$\Delta_h f  = \tau_h f -f $.

For functions $f \in \mathcal{S}(\Rn)$ the Fourier transform $\mathcal{F}$ and its inverse $\mathcal{F}^{-1} $ are defined via 
\[ \mathcal{F}(f)(\xi) = \il e^{-i x \cdot \xi} f(x) \, dx \,, \qquad
\mathcal{F}^{-1}(f)(x) = \il e^{i x \cdot \xi} f(\xi) \, \dd \xi \, ,\] 
where $\dd \xi=(2\pi)^{-n} d\xi$. When there is
ambiguity we use subscripts to indicate the variables with respect to which
the Fourier transform is taken, i.e., $\mathcal{F}(f)$ would be written as
$\mathcal{F}_{x\mapsto \xi }(f)$. Finally, $\mathcal{F}\colon \SD'(\Rn)\to
\SD'(\Rn)$ is defined by duality and $D_{x_j}:=\frac1{i}\partial_{x_j}$,
$j=1,\ldots, n$, where $\partial_{x_j}$ is the usual partial derivative. $D_x$ denotes the vector $(D_{x_1}, \ldots, D_{x_n})$.

We use a dyadic partition of unity $\varphi_j\in C_0^\infty(\Rn)$, $j\in\N_0$, which satisfies $\supp \varphi_0 \subset B_2(0)$ and $\supp \varphi_j \subset \{2^{j-1} \leq |\xi|\leq 2^{j+1}\}$ for $j\in \N$.  Then
the H\"older-Zygmund space $\mathcal{C}^s(\Rn)$, $s>0$, consists of all $f\in \SD'(\Rn)$ satisfying
\begin{equation*}
  \|f\|_{C^s} = \sup\{ 2^{ks}\|\varphi_k(D_x)f\|_{L^\infty}: k \in\N_0\} < \infty,
\end{equation*}
where 
\begin{equation*}
  \varphi_k(D_x) f= \F^{-1}\left[\varphi_k(\xi) \F[f](\xi)\right].
\end{equation*}
Note that $\mathcal{C}^s(\Rn)= B^s_{\infty \infty}(\Rn)$, where $B^s_{pq}(\Rn)$, $s\in \R$, $1\leq p,q\leq \infty$, denotes the usual Besov space. Moreover, it is well-known that $\sC^s(\Rn)=C^s(\Rn)$ for $s\in \R_+ \setminus \N $, cf. \cite[Appendix~A]{TaylorNonlinearPDE} or Triebel~\cite[Section~2.7]{Triebel1}.

The closure of $C_0^\infty(\Rn)$ in $\mathcal{C}^s(\Rn)$ is denoted by $\mathcal{C}^s_0(\Rn)$. 
We will use the following sufficient criterion for a function to belong to $\mathcal{C}^s_0(\Rn)$:
\begin{prop}\label{prop:Criterion}
  Let $0<s<s'<1$. Then every $f\in \sC^{s'}(\Rn)$ satisfying
  \begin{equation}
    \label{eq:Decay}
    \lim_{R\to\infty}\|f\|_{C^s(\Rn\setminus B_R(0))} =0
  \end{equation}
  belongs to $\mathcal{C}^s_0(\Rn)$.
\end{prop}
\begin{proof}
  Let $\varphi_\eps(x) = \eps^{-n}\varphi(\eps^{-1}x)$, $\varphi\in
  C_0^\infty(\Rn)$ with $\int\varphi(x) dx=1,$ be a standard mollifier. Then
  $\varphi_\eps \ast f\to_{\eps\to 0} f$ in $\sC^s(\Rn)$ since $f\in
  \sC^{s'}(\Rn)$. Moreover, (\ref{eq:Decay}) implies that each
  $\varphi_\eps\ast f$ can be approximated by smooth, compactly supported functions up to an arbitrarily small error in $C^s(\Rn)$. This proves the proposition.
\end{proof}

\subsection{Pseudodifferential Operators with Non-Smooth Symbols}

In the following, the principal part of the L\'evy-type operator will be represented as pseudodifferential operator with a symbol of the following kind:
\begin{defn}
  Let $n,n'\in \N$, $N\in\N_0$, $m\in \R$, and let $\tau \in (0,1)$. Then a function $p\colon \R^{n'}\times \Rn \to \C$ belongs to $C^\tau S^m_{1,0;N}(\R^{n'};\Rn)$ if $p(x,\xi)$ is H\"older continuous w.r.t. $x\in \R^{n'}$, $N$-times continuously differentiable w.r.t. $\xi \in \Rn$ and satisfies
  \begin{equation}
    \label{eq:SymbolEstim1}
    \|\partial_\xi^\beta p(.,\xi)\|_{C^\tau(\Rn)} \leq C \weight{\xi}^{\alpha-|\beta|}
  \end{equation}
  uniformly in $\xi\in \Rn$ and for all $|\beta|\leq N$. Moreover, let
  \begin{equation*}
    \|p\|_{C^\tau S^m_{1,0;N}}:= \sup_{\xi\in \Rn, |\beta|\leq N } \weight{\xi}^{-\alpha+|\beta|}\|\partial_\xi^\beta p(.,\xi)\|_{C^\tau(\Rn)}.
  \end{equation*}
\end{defn}
\begin{rem}
Note that $\bigcap\limits_{\tau>0, N \in \N} C^\tau S^m_{1,0;N}(\R^{n};\Rn)$ coincides with the classical symbol class $S^m_{1,0}(\R^{n};\R^n)$ as defined in \cite{Kum81}. A first treatment of pseudodifferential symbols which are merely H\"older continuous in the space variable $x$ and the associated operators was done by Kumano-Go and Nagase~\cite{KumanoGoNagase}. Further results and many references can be found in the monographs by Taylor~\cite{TaylorNonlinearPDE,ToolsForPDE}.  
\end{rem}

For $a=a(x,y,\xi)\in C^\tau S^m_{1,0;N}(\Rn\times \Rn;\Rn)$ we define the associated \emph{pseudodifferential operator in $(x,y)$-form} (formally) by
\begin{equation}\label{eq:DefnPsDO}
  a(x,D_x,x) f := \int_{\Rn}\int_{\Rn} e^{i(x-y)\cdot \xi} a(x,y,\xi) f(y) \sd y \dd \xi.
\end{equation}
So far, it is not clear whether $a(x,D_x,x) f$ in (\ref{eq:DefnPsDO}) is well-defined even for $f\in C^\infty_0(\Rn)$. This will be clarified later in each particular situation we have to deal with.

\begin{rem}
In order to underline the connection between the operator $a(x,D_x,x)$ and the corresponding symbol $a(x,y,\xi)$ we write $a(x,\xi,y)$ instead of $a(x,y,\xi)$ in the sequel. 
\end{rem}

In the special case that $a(x,\xi,y)=p(x,\xi)$, $p\in S^m_{1,0;N}(\Rn;\Rn)$, and $f\in \SD(\Rn)$,  the operator in (\ref{eq:DefnPsDO}) is well-defined as iterated integrals and coincides with
\begin{equation*}
  p(x,D_x) f= \int_{\Rn} e^{ix\cdot \xi} p(x,\xi) \hat{f}(\xi) \dd \xi,
\end{equation*}
which is a pseudodifferential operator in \emph{$x$-form}. The adjoints of $x$-form pseudodifferential operators are the pseudodifferential operators in $y$-form, which corresponds to the case $a(x,\xi,y)=p(y,\xi)$, $p\in S^m_{1,0;N}(\Rn;\Rn)$, and is (formally) given by
\begin{equation*}
  p(D_x,x) f := \F^{-1}\left[\int_{\Rn} e^{-iy\cdot \xi} p(y,\xi) f(y) \sd y\right].
\end{equation*}
If $f\in \SD(\Rn)$, the inner integral defines is a bounded continuous function in $\xi \in\Rn$ and $p(D_x,x)$ is a well-defined operator $p(D_x,x)\colon \SD(\Rn) \to \SD'(\Rn)$.

\begin{rem}
  Working with non-smooth symbols it is important to distinguish between
  pseudodifferential operators in $x$-form and in $y$-form since the mapping
  properties are different, cf. Theorem~\ref{thm:MappingProperties} below. The principal part of the operator $L$ will be a pseudodifferential operator
  in $x$-form; but it is important to take the approximate resolvent $Q_\lambda=q_\lambda(D_x,x)\approx (\lambda-L)^{-1}$ as an operator in $y$-form, not in $x$-form. Otherwise the mapping properties of $Q_\lambda$ would not fit to $(\lambda-L)^{-1}\colon C^s_0(\Rn)\to C^{s+\alpha}_0(\Rn)$ for $0<s<\tau$. 
This technique was already successfully applied to the resolvent equation of the Stokes operator in suitable domains with non-smooth boundary, cf. \cite{HInftyInLayer, BIPUnboundedDomains}. 
An alternative way for a parametrix construction is described in \cite[Section
  6]{NonsmoothGreen}, where the operator is first reduced to a zero order operator and then the parametrix is constructed in $x$-form. The latter article deals with pseudodifferential \emph{boundary value problems}; but the construction also applies to pseudodifferential equations on $\Rn$.
\end{rem}

Mapping properties of pseudodifferential operators with non-smooth
coefficients have been studied by several authors starting  with the pioneering work of Kumano-Go and Nagase~\cite{KumanoGoNagase}, cf. Taylor~\cite{TaylorNonlinearPDE,ToolsForPDE} and the references given there. For our purposes we will use the following theorem, which is a consequence of the results by Marschall~\cite{Marschall2}.
\begin{thm}\label{thm:MappingProperties}
  Let $N>\frac{n}2$, $\tau \in (0,1)$, and let $p\in C^\tau S^m_{1,0;N}(\Rn;\Rn)$. Then
  \begin{equation}\label{eq:MappingXForm}
    p(x,D_x)\colon \sC^{s+m}_0 (\Rn) \to \sC^{s} (\Rn)\qquad \text{if}\ 0<s<\tau, s+m>0
  \end{equation}
   and
  \begin{equation}\label{eq:MappingYForm}
    p(D_x,x)\colon \sC^{s+m}_0 (\Rn) \to \sC^{s} (\Rn)\qquad \text{if}\ s> 0, 0<s+m<\tau
  \end{equation}
  are bounded operators. Moreover, the operator norms can be estimated by $C\|p\|_{C^\tau S^m_{1,0;N}}$, where $C$ is independent of $p\in C^\tau S^m_{1,0;N}(\Rn;\Rn)$.
\end{thm}
\begin{rem}
  Note that for an operator $p(x,D_x)$ in $x$-form the order of the range space $\sC^{s}$ is limited by the smoothness of the symbol in $x$. For the corresponding operator in $y$-form, $p(D_x,x)$, the order of the domain $\sC_0^{s+m}$ is limited by $\tau$.
\end{rem}
\begin{proof*}{of Theorem~\ref{thm:MappingProperties}}
  First of all, we note that the symbol class $C^\tau S^m_{1,0;N}(\Rn;\Rn)$ coincides with the symbol class $S^m_{1,0}(\tau,N)$ defined in \cite{Marschall2}.
  Moreover, if $f\in \SD(\Rn)$, then $p(x,D_x)f$ defined as above coincides with the definition in \cite{Marschall2} as a limit of operators obtained from a symbol decomposition, cf. proof of \cite[Proposition~2.4]{Marschall2}.
  Hence \cite[Proposition~2.4]{Marschall2} implies that 
  \begin{equation*}
    \|p(x,D_x)f\|_{\sC^s(\Rn)}\leq C \|f\|_{\sC^{s+m}(\Rn)}
  \end{equation*}
  for $f\in \SD(\Rn)$ provided that $0<s<\tau$ and $s+m>0$. 

By our definition of $p(D_x,x)\colon \SD(\Rn)\to \SD'(\Rn)$
\begin{equation*}
  \weight{p(D_x,x) f, g} = \int_{\Rn}\int_{\Rn} e^{ix\cdot \xi} p(y,\xi) f(y) \sd y \hat{g}(-\xi) \dd \xi = \int_{\Rn} f(x) q(x,D_x) g \sd x
\end{equation*}
for all $f,g\in \SD(\Rn)$
with $q(x,\xi) = p(x,-\xi)$. Because of \cite[Proposition~4.3]{Marschall2}, $q(x,D_x)^\ast \colon \sC^{s+m}(\Rn)\to \sC^s(\Rn)$ provided that $0<s+m<\tau$ and $s>0$. 

Finally, it is easy to observe that all estimates done in the proof of  \cite[Proposition~4.3]{Marschall2} are uniform for all $p\in C^\tau S^m_{1,0;N}(\Rn;\Rn)$ with $\|p\|_{C^\tau S^m_{1,0;N}}\leq 1$, which is nothing but the boundedness of the linear mapping from the symbol space $C^\tau S^m_{1,0;N}(\Rn;\Rn)$ into the corresponding space of linear operators.
\end{proof*}

The next important ingredient are kernel estimates of the Schwartz kernel associated to a pseudodifferential operator. We follow the presentation given in \cite[Chapter 6, Paragraph 4]{Ste93}. Given $a\in C^\tau S^m_{1,0;N}(\Rn\times \Rn;\Rn)$ we define for $j\in \N_0$
\begin{equation*}
  k_j(x,y,z):= \F^{-1}_{\xi\mapsto z}[a_j(x,.,y)], \qquad a_j(x,\xi,y):=a(x,\xi,y)\varphi_j(\xi),
\end{equation*}
where $\varphi_j$ is the Dyadic partition of unity introduced above.

First of all, we have
\begin{lem}\label{lem:EstimKj}
  Let $a\in C^\tau S_{1,0;N}^m(\Rn\times\Rn)$, $m\in\R$, $N\in\N_0$, $\tau\in (0,1)$, and let $k_j(x,y,z)$ be defined
  as above. Then
  \begin{equation}\label{eq:EstimKj}
    \|\partial_z^\alpha k_j(.,.,z)\|_{C^\tau (\Rn\times\Rn)}\leq C_{\alpha,M}\|a\|_{C^\tau S_{1,0;N}^m}
     |z|^{-M} 2^{j(n+m-M+|\alpha|)}
  \end{equation}
  for all $\alpha\in\N_0^n$, $M=0,\ldots,N$, where $C_{\alpha,M}$
  does not depend on $j\in\N_0$ and $a\in C^\tau S_{1,0;N}^m(\Rn\times\Rn;\Rn)$.
\end{lem}
\begin{proof}
  We start with 
  \begin{equation*}
    z^\gamma D_z^\alpha k_j(x,z) = \int_{\Rn} e^{ix\cdot \xi}D_\xi^\gamma
    [\xi^\alpha a_j(x,y,\xi)] \dd \xi
  \end{equation*}
  for all $\alpha,\gamma\in\N_0^n$. We estimate the integral on the right hand side from above. Firstly, the integrand is supported in
  the ball $\{|\xi|\leq 2^{j+1}\}$, which has volume bounded by a multiple of
  $2^{nj}$. Secondly, since the support is also limited by the condition $2^{j-1}\leq |\xi|$
  (when $j\neq 0$) and $c2^{j}\leq \weight{\xi}\leq C 2^{j}$ on $\{2^{j-1}\leq
  |\xi|\leq 2^{j+1}\}$,
  \begin{equation*}
    \left|D_\xi^\gamma
    [\xi^\alpha a_j(x,y,\xi)]\right|\leq
    C_{\alpha,\gamma} \|a\|_{C^\tau S_{1,0;N}^m}2^{j(m+|\alpha|-|\gamma|)}
  \end{equation*}
due to the symbol estimates of $\xi^\alpha a_j(x,y,\xi)\in
C^\tau S^{m+|\alpha|}_{1,0;N}(\Rn\times\Rn;\Rn)$. Hence
  \begin{equation*}
    |z^\gamma D_z^\alpha k_j(x,y,z)| \leq C_{\alpha,\gamma}\|a\|_{C^\tau S_{1,0;N}^m}
     2^{j(n+m+|\alpha|-M)}, \qquad \text{whenever}\ |\gamma|=M.
  \end{equation*}
  Taking the supremum over all $\gamma$ with $|\gamma|=M$, gives
  (\ref{eq:EstimKj}) with $C^\tau(\Rn\times\Rn)$ replaced by $C^0(\Rn\times \Rn)$. In order to get the same for $C^\tau(\Rn\times\Rn)$ one simply replaces all terms above by suitable differences.
\end{proof} 

Using the latter lemma, we are able to prove the following kernel estimate:
\begin{thm}\label{thm:kernelEstim}
  Let $a\in C^\tau S^m_{1,0;N}(\Rn\times\Rn;\Rn)$, $\tau\in (0,1)$, $m>-n$, and $N\in\N_0$ such that $N>n+m$ and let $k_j$ be defined as above. 
  Then for every $x,y,z\in \Rn, z\neq 0,$ 
  \begin{equation*}
    k(x,y,z):=  \sum_{j=0}^\infty k_j(x,y,z)
  \end{equation*}
  exists, converges uniformly in $x,y\in\Rn$, $|z|\geq \eps>0$, and satisfies
  \begin{equation*}
    \|\partial_z^\alpha k(.,.,z)\|_{C^\tau(\Rn\times\Rn)} \leq 
    \begin{cases}
      C_\alpha\|a\|_{C^\tau S^m_{1,0;N}} |z|^{-n-m-|\alpha|} & \text{for}\ |z|\leq 1 \\
C_\alpha\|a\|_{C^\tau S^m_{1,0;N}} |z|^{-N} & \text{for}\ |z|\geq 1
    \end{cases}
  \end{equation*}
  uniformly in $z\neq 0$
  for all $\alpha\in \N_0$ with $|\alpha| < N-n-m$,
  where $C$ is independent of $a\in C^\tau S^m_{1,0;N}(\Rn\times\Rn;\Rn)$.
\end{thm}
\begin{proof}
  First we consider the case when $0<|z|\leq 1$. We brake the  above sum
into two parts: the first where $2^j\leq |z|^{-1}$, the second where
$2^j>|z|^{-1}$. In order to estimate the first sum we use (\ref{eq:EstimKj}) 
with $M=0$:
\begin{equation*}
  \sum_{2^j\leq |z|^{-1}} \|\partial_z^\alpha k_j(.,.,z)\|_{C^\tau(\Rn\times \Rn)} \leq 
  C \|a\|_{C^\tau S^m_{1,0;N}}\sum_{2^j\leq |z|^{-1}} 2^{j(n+m+|\alpha|)},
\end{equation*}
where
\begin{equation*}
  \sum_{2^j\leq |z|^{-1}} 2^{j(n+m+|\alpha|)}= 
    O(|z|^{-n-m-|\alpha|}) 
\end{equation*} 
since $n+m+|\alpha|>0$.

Next, for the second sum, we use again (\ref{eq:EstimKj}) with $M=N$ and get the estimate
\begin{eqnarray*}
  \sum_{2^j>|z|^{-1}} \|\partial_z^\alpha k_j(.,z)\|_{C^\tau(\Rn\times \Rn)} &\leq&
  C_{\alpha} \|a\|_{C^\tau S^m_{1,0;N}}  |z|^{-M}\sum_{2^j>|z|^{-1}} 2^{j(n+m+|\alpha|-M)}\\ 
  &\leq& 
  C'_{\alpha}\|a\|_{C^\tau S^m_{1,0;N}} |z|^{-n-m-|\alpha|}.
\end{eqnarray*}
Finally, we consider the situation $|z|\geq 1$. Since $N>n+m+|\alpha|$, 
(\ref{eq:EstimKj}) shows that 
  \begin{eqnarray*}
  \sum_{j=0}^\infty \|\partial_z^\alpha k_j(.,z)\|_{C^\tau(\Rn\times \Rn)}&\leq& C_{\alpha} |z|^{-N}\|a\|_{C^\tau S^m_{1,0;N}} \sum_{j=0}^\infty 2^{j(n+m-N+|\alpha|)}\\ &\leq& C'_{\alpha}\|a\|_{C^\tau S^m_{1,0;N}} |z|^{-N}. 
\end{eqnarray*}
Hence the proof is complete.
\end{proof}

The following corollary shows that (\ref{eq:MappingXForm}) can be improved to  $p(x,D_x)\colon \sC^{s+m}_0 (\Rn) \to \sC^{s}_0 (\Rn)$ under the same assumptions.

\begin{cor}\label{cor:Mapping}
  Let $N>n+m$, $\tau \in (0,1)$, let $p\in C^\tau S^m_{1,0;N}(\Rn;\Rn)$, and let $f\in C_0^\infty(\Rn)$. Then $p(x,D_x)f\in \lh^s(\Rn)$ for all $0<s<\tau$ with $s+m>0$ and
  $p(D_x,x) f\in \lh^s(\Rn)$ provided that $0<s+m<\tau$ and $s>0$.
\end{cor}
\begin{proof}
  For simplicity we only treat the case of the operator in $x$-form. The other case is treated in the same way.

  Fix $0<s<\tau$ with $s+m>0$ and choose $s'\in (s,\tau)$. Then $p(x,D_x)f\in \sC^{s'}(\Rn)$ due to Theorem~\ref{thm:MappingProperties}.
  Hence, using Proposition~\ref{prop:Criterion}, it is sufficient to show (\ref{eq:Decay}). Because of Theorem~\ref{thm:kernelEstim} with $a(x,\xi,y)=p(x,\xi)$,
  \begin{eqnarray*}
    p(x,D_x) f&=& \sum_{j=0}^\infty p_j(x,D_x) f 
    = \sum_{j=0}^\infty \int_{\Rn} k_j(x,x-y) f(y) \sd y\\
 &=& \int_{\Rn} k(x,x-y) f(y) \sd y\qquad \text{for all}\ x\not\in \supp f.
  \end{eqnarray*}
  Using the kernel estimate stated in Theorem~\ref{thm:kernelEstim}, one easily verifies (\ref{eq:Decay}).  
\end{proof}


Recall that, if $a \in S^m_{1,0}(\Rn\times \Rn;\Rn)$ is a smooth symbol, then by the results of the classical theory of pseudodifferential operators
\[ a(x,D_x,x) = p(x,D_x) \,, \]
where $p \in S^m_{1,0}(\Rn\times \Rn;\Rn)$ and
\[ p(x,\xi) = a(x,\xi,x) + r(x,\xi)  \,, \]
with $r \in S^{m-1}_{1,0}(\Rn;\Rn)$, see \cite[Chapter 2, Section 3]{Kum81}. In the case $a \in C^\tau S^m_{1,0;N}(\Rn\times \Rn;\Rn),0 \leq \tau \leq m$, the following result can be applied to
\[ r(x,\xi,y)=a(x,\xi,y) - a(x,\xi,x) \,.\] 

\begin{prop}\label{prop:RedOp}
  Let $r\in C^\tau S^m_{1,0;N}(\Rn\times \Rn;\Rn)$, where $\tau \in (0,1)$,
  $0\leq m < \tau$, and $N=n+1$. Moreover, we assume that $r(x,\xi,x)=0$. Then
  \begin{equation*}
    r(x,D_x,x) := \sum_{j=0}^\infty r_j(x,D_x,x) 
  \end{equation*}
  converges absolutely in $\mathcal{L}(\sC^s(\Rn))$ for each $0<s<\tau-m$ and
  satisfies
  \begin{equation}\label{eq:EstimRedOP}
    \|r(x,D_x,x)\|_{\mathcal{L}(\sC^s_0(\Rn))}\leq C \|r\|_{C^\tau S^m_{1,0;N}},
  \end{equation}
  where $C$ does not depend on $r\in C^\tau S^m_{1,0;N}(\Rn\times \Rn;\Rn)$.
  Moreover, $r(x,D_x,x)$ maps $\lh^s(\Rn)$ into itself.
\end{prop}
\begin{proof}
  First we denote 
  \begin{equation*}
    r^M(x,D_x,x) f:= \sum_{j=0}^M r_j(x,D_x,x) f.
  \end{equation*}
  Using that
  \begin{equation*}
    r_j(x,D_x,x) = \int_{\Rn}\int_{\Rn} k_j(x,y,x-y) f(y) \sd y, \quad f\in \SD(\Rn),
  \end{equation*}
  we have
  \begin{equation*}
    r^M(x,D_x,x) = \int_{\Rn} k^M(x,y,x-y) f(y) \sd y, \quad f\in \SD(\Rn),
  \end{equation*}
  with $k^M(x,y,z):= \sum_{j=0}^M k_j(x,y,z)$. Note that $k^M(x,x,z)=k_j(x,x,z)=0$ since $r(x,\xi,x)=0$. By the proof of Theorem~\ref{thm:kernelEstim} it is obvious that 
  \begin{equation*}
    \| k^M(.,z)\|_{C^\tau(\Rn\times\Rn)} \leq 
    \begin{cases}
      C\|r\|_{C^\tau S^m_{1,0;N}} |z|^{-n-m} & \text{if}\ |z|\leq 1\\
C\|r\|_{C^\tau S^m_{1,0;N}} |z|^{-n-1} & \text{if}\ |z|\geq 1,
    \end{cases}
  \end{equation*}
  uniformly in $z\neq 0$ and $M\in\N$. But this implies
  \begin{eqnarray}\nonumber
    |k^M(x,y,x-y)|&=&|k^M(x,y,x-y)-k^M(x,x,x-y)|\\\label{eq:EstimkM}
    &\leq& C \|r\|_{C^\tau S^m_{1,0;N}}|x-y|^{-n-m+\tau}(1+|x-y|)^{m-1}.
  \end{eqnarray}
  Hence Lebesgue's theorem on dominated convergence implies that
  \begin{equation*}
    r(x,D_x,x) f = \lim_{M\to \infty} r^M(x,D_x,x) f = \int_{\Rn} k(x,y,x-y) f(y) \sd y
  \end{equation*}
  exists 
  for every $x\in \Rn$ and $f\in L^\infty(\Rn)$. Moreover, since (\ref{eq:EstimkM}) holds for $k(x,y,x-y)$ as well, we conclude
  \begin{equation}\label{eq:LinftyEstim}
    \|r(x,D_x,x) \|_{\mathcal{L}(L^\infty(\Rn))}\leq C\|r\|_{C^\tau S^m_{1,0;N}}. 
  \end{equation}
  In order to prove (\ref{eq:EstimRedOP}), we use the relation
  \begin{equation*}
    \Delta_h r(x,D_x,x)f = r(x,D_x,x) (\Delta_h f) + \int_{\Rn} k_h(x,y,x-y) f(y+h) \sd y,  
  \end{equation*}
  where $(\Delta_h f)(x) = f(x+h)-f(x)$, $h\in\Rn$, and
  \begin{equation*}
    k_h(x,y,z)= k(x+h,y+h,z)-k(x,y,z).
  \end{equation*}
  Moreover, $k_h(x,y,z)$ is the kernel belonging to $r_h(x,D_x,x)$ with $r_h(x,y,\xi)= r(x+h,\xi,y+h)- r(x,\xi,y)$ and it is easy to prove that 
  \begin{equation*}
    \|r_h\|_{C^{\tau-s}S^m_{1,0;N}}\leq C|h|^s \|r_h\|_{C^{\tau}S^m_{1,0;N}}
  \end{equation*}
  uniformly in $h\in \Rn$ for each $0<s<\tau$.
  Hence using (\ref{eq:LinftyEstim}) for $r$ and $r_h$, we conclude that
  \begin{eqnarray*}
    \|\Delta_h r(x,D_x,x)f\|_{L^\infty} &\leq& C \|r\|_{C^\tau S^m_{1,0;N}} \|\Delta_h f\|_{L^\infty} + C \|r_h\|_{C^{\tau-s} S^m_{1,0;N}} \|f\|_{L^\infty}\\
    &\leq& C \|r\|_{C^\tau S^m_{1,0;N}} \|f\|_{C^s(\Rn)} |h|^s
  \end{eqnarray*}
  for $0<s<\tau-m$. This finishes the proof of (\ref{eq:EstimRedOP}). The last statement is proved in the same way as in Corollary~\ref{cor:Mapping}.
\end{proof}

\subsection{Application to the Resolvent Equations}

In this section we construct an approximate resolvent $Q_\lambda$ to a L\'evy-type operator $L$ as introduced in (\ref{eq:DefnL}),(\ref{eq:DefnL2}). Here $Q_\lambda=q_\lambda(D_x,x)$ is a pseudodifferential operator obtained by inverting the symbol of the principal part of $\lambda-L$. 

More precisely, because of the assumption on the kernel, we have a decomposition
\begin{equation*}
  L u(x) = L^1 u(x) + L^2 u(x),\qquad u\in \SD(\Rn),
\end{equation*}
where $L^j$ denotes the same kind of operator with kernel $k^j$, $j=1,2$.
Here $L^1$ can be considered as principle part and $L^2$ is of lower order in the following sense:
\begin{lem}
  Let $L^2$ be as above. Then $L^2$ extends to a bounded operator
  $L^2\colon \lh^{s+\alpha''}(\Rn)\to \lh^s (\Rn)$ for any $\alpha''>\alpha'$ and $0<s<\tau$ provided that $s+\alpha''> 1$ if $\alpha\geq 1$.
\end{lem}
\begin{proof}
  First of all, if $u\in \sC^{s'}(\Rn)$ and $1<s'< 2$, then
  \begin{equation}\label{eq:EstimUx}
    \left|u(x+y)-u(x)-y\cdot\nabla u(x)
  \right| \leq C\|u\|_{\sC^{s'}(\Rn)} |y|^{s'},\quad |y|\leq 1. 
  \end{equation}
  First we assume that $1\leq \alpha'<\alpha <2$. 
  Then (\ref{eq:EstimUx}) with $s'=\alpha''$ yields
  \begin{eqnarray}\nonumber
    \lefteqn{\|L^2u\|_{L^\infty(\Rn)}} \\\label{eq:L2Estim1} 
    &\leq& C
    \left(\sup_{x\in\Rn,|y|\leq 1} |y|^{n+\alpha'}|k_2(x,y)| + \int_{|y|\geq 1} \|k_2(.,y)\|_\infty \sd y \right)
\|u\|_{\sC^{\alpha''}(\Rn)}
  \end{eqnarray}
  with a constant $C$ independent of $k_2$. Moreover, 
  \begin{equation}\label{eq:Diff}
    \Delta_h (L^2 u) = L^2 (\Delta_h u) + L_h^2 (\tau_h u),
  \end{equation}
  where $L_h^2$ is the L\'evy-type operator with kernel
  $
    k^2_h(x,y):= k^2(x+h,y)-k^2(x,y)
  $.
  By the assumptions on the kernel, 
  \begin{equation*}
    \sup_{x\in\Rn,|y|\leq 1} |y|^{n+\alpha'}|k^2_h(x,y)| + \int_{|y|\geq 1} \|k^2_h(.,y)\|_\infty \sd y \leq C|h|^s 
  \end{equation*}
  uniformly in $h\in\Rn$. Therefore using
(\ref{eq:L2Estim1}) with $L^2$ replaced by holds for $L^2_h$ and $k_2$ replaced by $k^2_h$ we conclude
\begin{equation*}
  \|L^2_h (\tau_h u)\|_{L^\infty(\Rn)}\leq C |h|^s\|u\|_{\sC^{\alpha''}(\Rn)}. 
\end{equation*}
Hence, using the inequality above, (\ref{eq:Diff}), and (\ref{eq:L2Estim1}), we conclude
\begin{equation*}
  \|\Delta_h (L^2 u)\|_{L^\infty(\Rn)}\leq C\left( \|\Delta_h u\|_{\sC^{\alpha''}(\Rn)}+ |h|^s\|u\|_{C^{\alpha''}(\Rn)}\right) \leq C h^s \|u\|_{\sC^{s+\alpha''}(\Rn)},
\end{equation*}
where we have used $\|\Delta_h u\|_{\sC^{\alpha''}(\Rn)}\leq C|h|^s \|
u\|_{\sC^{s+\alpha''}(\Rn)}$. The latter inequality can be easily proved by first proving the
cases $s=0,1$ and then using interpolation. Hence $L^2\colon \sC^{s+\alpha''}(\Rn)\to \sC^s (\Rn) $. 

Secondly, if $0<\alpha<1$, then the proof above is easily modified using
  \begin{equation*}
    \left|u(x+y)-u(x)\right| \leq C\|u\|_{\sC^{s'}(\Rn)} |y|^{s'},\quad |y|\leq 1, 
  \end{equation*}
  for $u\in \sC^{s'}(\Rn)$ and $s'\in (0,1)$ instead of (\ref{eq:EstimUx}).

It remains to consider the case $0\leq \alpha'<1\leq \alpha$. Using (\ref{eq:EstimUx}) with $s'=s+\alpha''\in (1,2)$ we conclude as before
  \begin{eqnarray}\nonumber
    \lefteqn{\|L^2u\|_{L^\infty(\Rn)}} \\\label{eq:L2Estim2} 
    &\leq& C
    \left(\sup_{x\in\Rn,|y|\leq 1} |y|^{n+\alpha'}|k_2(x,y)| + \int_{|y|\geq 1} \|k_2(.,y)\|_\infty \sd y \right)
\|u\|_{\sC^{s+\alpha''}(\Rn)}
  \end{eqnarray}
  with a constant $C$ independent of $k_2$. We use again (\ref{eq:Diff}). The second term can be estimated in the same manner as before to obtain
\begin{equation*}
  \|L^2_h (\tau_h u)\|_{L^\infty(\Rn)}\leq C |h|^s\|u\|_{\sC^{s+\alpha''}(\Rn)}. 
\end{equation*}
But the first term in (\ref{eq:Diff}) has to be estimated differently: Using (\ref{eq:EstimUx}) with $u$ replaced by $\Delta_h u$, we have on one hand
\begin{eqnarray*}
  \lefteqn{\left|\Delta_h u(x+y)-\Delta_h u(x)-y\cdot\nabla \Delta_h u(x)\right|}\\ 
  &\leq& C\|\Delta_h u\|_{\sC^{s+\alpha''}(\Rn)} |y|^{s+\alpha''}
  \leq C'\|u\|_{\sC^{s+\alpha''}(\Rn)} |y|^{s+\alpha''},\quad |y|\leq 1. 
\end{eqnarray*}
On the other hand
\begin{eqnarray*}
  \lefteqn{\left|\Delta_h u(x+y)-\Delta_h u(x)-y\cdot\nabla \Delta_h u(x)
      \right|}\\
  &\leq& C\|\Delta_h u\|_{C^1(\Rn)} |y|\leq C'|y||h|^{s+\alpha''-1}\|u\|_{\sC^{s+\alpha''}(\Rn)},\quad |y|,|h|\leq 1. 
\end{eqnarray*}
Interpolation of both inequalities yields
\begin{equation*}
  \left|\Delta_h u(x+y)-\Delta_h u(x)-y\cdot \nabla \Delta_h u(x)
    \right|
  \leq C |h|^s|y|^{\alpha''}\|u\|_{\sC^{s+\alpha''}(\Rn)}\quad 
\end{equation*}
uniformly in $|h|,|y|\leq 1$.
With this inequality
\begin{equation*}
  \|L^2 \Delta_h u \|_{L^\infty(\Rn)}\leq C |h|^s\|u\|_{\sC^{s+\alpha''}(\Rn)}, \quad |h|\leq 1,
\end{equation*}
is proved in the same way as before.

Finally, if $f\in C_0^\infty(\Rn)$, one easily proves $L^2 f\in \sC^s_0 (\Rn)$ with the aid of Proposition~\ref{prop:Criterion} and (\ref{eq:kernelEstim5}).
\end{proof}

For the principal part $L^1$, we use
\begin{alignat*}{1}
  u(x+y)-u(x)-y\cdot \nabla u(x)
   &= \F^{-1}_{\xi\mapsto x}\left[ \left(e^{iy\cdot \xi}-1-i\xi\cdot
  y\right) \hat{u}(\xi)\right],\\
  u(x+y)-u(x) &= \F^{-1}_{\xi\mapsto x}\left[ \left(e^{iy\cdot \xi}-1\right) \hat{u}(\xi)\right].
\end{alignat*}
Hence 
$L^1$ can be represented as a pseudodifferential operator
\begin{eqnarray*}
  L^1u(x) &=& 
  \int_{\Rn} e^{ix\cdot\xi} p(x,\xi) \hat{u}(\xi) \dd \xi,
\end{eqnarray*}
where
\begin{alignat*}{2}
  p(x,\xi)&:= \int_{\R^n} \left(e^{iy\cdot \xi}-1-i\xi\cdot
  y\right)k_1(x,y)\sd y &\quad & \text{if} \ \alpha\in [1,2),\\
  p(x,\xi)&:= \int_{\R^n} \left(e^{iy\cdot \xi}-1\right)k_1(x,y)\sd y &\ & \text{if} \ \alpha\in (0,1).
\end{alignat*}

The following lemma shows that $p$ is a symbol in the class studied above.
\begin{lem}\label{lem:kernelEstim}
  Let $k_1\colon\R^n \times \R^n \to \R$ be $N$-times differentiable w.r.t the
  second variable satisfying 
  \begin{equation}
    \|\partial_y^\beta k_1(.,y)\|_{C^\tau(\Rn)}\leq C |y|^{-n-\alpha-|\beta|}
  \end{equation}
  for all $0<|y|\leq 2$ and $|\beta|\leq N$ and $k(x,y)=0$ for $|y|\geq 2$.  Then $p\in C^\tau S^\alpha_{1,0;N}(\Rn;\Rn)$ where $p$ is defined as above.
\end{lem}
\begin{proof}
  We denote $f(s)= e^{is}-1-is$, $s\in \R$, if $\alpha\in [1,2)$ and $f(s)= e^{is}-1$, $s\in \R$, if $\alpha \in (0,1)$.
  Let $\gamma,\beta\in\N_0^n$ with $m=|\gamma|=|\beta|\leq N$. Then
  \begin{equation*}
    \partial_{\xi}^\beta\left(\xi^\gamma f(y\cdot\xi)\right) = 
\partial_{\xi}^\beta\left(\partial_y^\gamma F^m(y\cdot\xi)\right) =
\partial_y^\gamma\left( \partial_{\xi}^\beta F^m(y\cdot\xi)\right) =
 \partial_y^\gamma\left( y^\beta f(y\cdot\xi)\right)
  \end{equation*}
  where $F^m$ denotes the $m$-th primitive of $f$. Therefore
  \begin{eqnarray*}
     \partial_{\xi}^\beta\left(\xi^\gamma p(x,\xi)\right)
&=& \int_{\Rn} \partial_y^\gamma \left( y^\beta f(y\cdot
      \xi)\right) k_1(x,y) \sd y \\
&=& (-1)^m\int_{\Rn} y^\beta f(y\cdot \xi) \partial_y^\gamma k_1(x,y) \sd y\\ 
&=& (-1)^m|\xi|^{-n-m}\int_{\Rn} z^\beta f\left(z\cdot \frac{i\xi}{|\xi|}\right) (\partial_y^\gamma k_1)\left(x,\frac{z}{|\xi|}\right) \sd z 
  \end{eqnarray*}
  Hence 
  \begin{eqnarray*}
    \left\|\partial_{\xi}^\beta \left(\xi^\gamma p(.,\xi)\right)\right\|_{C^\tau(\Rn)}
    &\leq & C |\xi|^{-n-m}\int_{\Rn} |z|^m \frac{|z|^j}{1+|z|^j}
    \left|\frac{z}{|\xi|}\right|^{-n-\alpha-m} \sd z \\
&\leq& C'|\xi|^{\alpha},
  \end{eqnarray*}
  where $j=2$ if $\alpha\geq 1$ and $j=1$ else.
  Since $\beta,\gamma\in \N_0^n$ with $|\beta|=|\gamma|\leq N$ are arbitrary, this implies
  \begin{equation*}
    \left\|\xi^\gamma \partial_{\xi}^\beta p(.,\xi)\right\|_{C^\tau(\Rn)}\leq C|\xi|^\alpha
  \end{equation*}
  for all $|\beta|=|\gamma|\leq N$, which is easy to prove by induction. Hence
  \begin{equation*}
    \left\|\partial_{\xi}^\beta p(.,\xi)\right\|_{C^\tau(\Rn)}\leq C|\xi|^{\alpha-|\beta|}
  \end{equation*}
  since $\gamma\in\N_0^n$ with $|\gamma|=|\beta|$ is arbitrary. 
\end{proof}

Hence (\ref{eq:kernelEstim1}), Lemma~\ref{lem:kernelEstim}, Theorem~\ref{thm:MappingProperties}, and Corollary~\ref{cor:Mapping} imply that
\begin{equation*}
  p(x,D_x) \colon \sC^{s+\alpha}_0(\Rn) \to \sC^{\alpha}_0(\Rn)
\end{equation*}
for all $0<s<\tau$. Moreover, (\ref{eq:kernelEstim2}) implies
\begin{equation*}
  -\re p(x,\xi) = \int_{\Rn} (1-\cos y\cdot \xi) k_1(x,y) \sd y
  \geq c\int_{B_2(0)} (1-\cos y\cdot \xi) |y|^{-n-\alpha} \sd y
  \geq C |\xi|^\alpha
\end{equation*}
for all $|\xi|\geq 1$ and $-\re p(x,\xi)\geq 0$ for all $\xi\in\Rn$.
Since $|p(x,\xi)|\leq C \weight{\xi}^\alpha$, we conclude that
\begin{equation*}
  \left|\frac{\im p(x,\xi)}{\re p(x,\xi)}\right|\leq M
\end{equation*}
uniformly in $|\xi|\geq 1$. Thus $p(x,\xi) \in \C\setminus \Sigma_\delta$ for $\delta:=\pi - \arctan M >\frac{\pi}2$ and for all $|\xi|\geq 1$.

Hence, we can define
\begin{equation*}
  q_\lambda(y,\xi) := (\lambda-p(y,\xi))^{-1},\qquad y,\xi\in\Rn,\lambda \in \Sigma_{\delta'}, |\lambda|\geq R,
\end{equation*}
for $0<\delta'<\delta$ and $R>\sup_{x\in\Rn, |\xi|\leq 1} |p(x,\xi)|$.

Since $p\in C^\tau S^{\alpha}_{1,0;N}(\Rn;\Rn)$, we have $q_\lambda\in C^\tau S^{-\alpha}_{1,0;N}(\Rn;\Rn)$. More precisely, the following lemma holds: 
\begin{lem}\label{lem:ParametrixSymbol}
  Let $q_\lambda$, $\delta$ be defined as above and $\lambda \in\Sigma_{\delta'}$ where $\delta'\in(0,\delta)$ is arbitrary.  Then there is some $R>0$ such that $q_\lambda \in C^\tau S^{-\alpha}_{1,0;N}$ for all $\lambda \in\Sigma_{\delta'}$ with $|\lambda|\geq R$. Moreover, for each $\alpha'\in [0,\alpha]$
  \begin{equation*}
    \|q_\lambda\|_{C^\tau S^{-\alpha'}_{1,0;N}} \leq C_{\delta'} (1+|\lambda|)^{-\frac{\alpha-\alpha'}\alpha}
  \end{equation*}
  uniformly in $\lambda \in\Sigma_{\delta'}$ with $|\lambda|\geq R$.
\end{lem}
\begin{proof}
  First of all, by a simple geometric observation
  \begin{equation*}
    |\lambda -z|\geq c_{\delta'} \max \{|\lambda|, |z|\}\quad
    \text{if}\ \lambda \in \Sigma_{\delta'}, z\in \C\setminus \Sigma_{\delta}
  \end{equation*}
  provided that $0<\delta'<\delta$. As seen above $p(x,\xi)\in \C\setminus\Sigma_\delta$
  for $|\xi|\geq 1$ and some $\delta>\frac{\pi}2$ and $|p(x,\xi)|\geq c|\xi|^\alpha$ for $|\xi|\geq 1$. Hence
  \begin{equation}\label{eq:EstimBelow1}
    |\lambda - p(x,\xi)|\geq c_{\delta'} \max \{|\lambda|, |\xi|^\alpha\}
  \end{equation}
  for all $|\xi|\geq 1$ and $\lambda\in \Sigma_{\delta'}$ with $0<\delta'<\delta$ arbitrary. Moreover, since $|p(x,\xi)|\leq C$ for all $|\xi|\leq 1$ and $x\in \Rn$, we conclude that (\ref{eq:EstimBelow1}) holds for all $\xi\in\Rn$ and $\lambda\in \Sigma_{\delta'}$ with $|\lambda|\geq R$ for some $R>0$ sufficiently large. 
Using this, $p\in C^\tau S^\alpha_{1,0;N}(\Rn;\Rn)$, and the chain rule, one derives in a straight-forward manner that
  \begin{equation*}
    \|\partial_\xi^\beta q_\lambda(.,\xi)\|_{C^\tau(\Rn)}\leq C_{\delta'} \frac{\weight{\xi}^{-|\beta|}}{|\lambda|+|\xi|^\alpha}\leq C_{\delta'}|\lambda|^{-\frac{\alpha-\alpha'}\alpha}\weight{\xi}^{-\alpha'-|\beta|} 
    \end{equation*}
  uniformly in $\xi\in\Rn$ and $\lambda\in \Sigma_{\delta'}$, $|\lambda|\geq
  R>0$ and for all $|\beta|\leq N$, which proves the statement.
\end{proof}

Application of Theorem~\ref{thm:MappingProperties}, Corollary~\ref{cor:Mapping} and the lemma above gives:
\begin{cor}\label{cor:MappingOfParametrix}
  Let $q_\lambda, \delta, \delta'$ be as above and let $0<s<\tau$. Then $q_\lambda(D_x,x)\colon \lh^{s}(\Rn) \to \lh^{s+\alpha}(\Rn)$ is a bounded linear operator, which satisfies
  \begin{equation*}
    \|q_\lambda(D_x,x)\|_{\mathcal{L}(\sC^{s}_0(\Rn), \sC^{s+\alpha'}_0(\Rn))}
    \leq C_{\delta'}|\lambda|^{-\frac{\alpha-\alpha'}\alpha}\qquad
      \text{for all} \ \lambda\in \Sigma_{\delta`}, |\lambda|\geq R,
  \end{equation*}
  for all $0\leq \alpha'\leq \alpha$ with some sufficiently large $R>0$.
\end{cor}

Now we are in the position to prove the following key lemma.
\begin{lem}\label{lem:Parametrix}
  Let $q_\lambda, \delta, \delta'$ be as above and let $0<s<\tau$. Then
  \begin{equation*}
    (\lambda-p(x,D_x))q_\lambda(D_x,x) = I - R_\lambda
  \end{equation*}
  with
  \begin{equation*}
    \|R_\lambda\|_{\mathcal{L}(\sC^s_0(\Rn))}\leq C_{\delta'}|\lambda|^{-\eps}
  \end{equation*}
  uniformly in $\lambda\in \Sigma_{\delta'}$ with $|\lambda|\geq M$ 
  for sufficiently large $M>0$ and some $\eps>0$ depending on $s,\tau$.
\end{lem}
\begin{proof}
First of all, for each $f\in C_0^\infty(\Rn)$, $q_\lambda(D_x,x) f\in \sC^{s'+\alpha}(\Rn)$ with $s<s'<\tau$. 
We conclude
\begin{equation*}
\sum_{j=0}^N \varphi_j(D_x) q_\lambda(D_x,x) f\to  q_\lambda (D_x,x) f \qquad \text{ in } \sC^{s+\alpha}(\Rn) \text{ as } {N\to\infty}.
\end{equation*}
 Therefore
\begin{equation*}
 q_\lambda(D_x,x) f=  \sum_{j=0}^\infty \varphi_j(D_x) q_\lambda(D_x,x)  f
 =\sum_{j=0}^\infty q_{\lambda,j}(D_x,x)f
\end{equation*}
where $q_{\lambda,j}(\xi,y)= q_\lambda(\xi,y)\varphi_j(\xi)$.
Hence
  \begin{eqnarray*}
    (\lambda-p(x,D_x))q_\lambda(D_x,x)f &=& \sum_{j=0}^\infty (\lambda - p(x,D_x)) q_{\lambda,j}(D_x,x)f\\
    &=& f+ \sum_{j=0}^\infty a_{\lambda,j}(x,D_x,x) f, 
  \end{eqnarray*}
  where $a_{\lambda,j}(x,y,\xi)= a_\lambda(x,\xi,y)\varphi_j(\xi)$ and
  \begin{equation*}
    a_{\lambda}(x,y,\xi) = \frac{\lambda-p(x,\xi)}{\lambda-p(y,\xi)} -1
    = (p(y,\xi)-p(x,\xi))q_\lambda(y,\xi).
  \end{equation*}
  Using Lemma~\ref{lem:ParametrixSymbol}, we conclude
  \begin{eqnarray*} 
    \|a_\lambda\|_{C^\tau S^{\alpha-\alpha'}_{1,0;N}}\leq C \|p\|_{C^\tau S^{\alpha}_{1,0;N}}\|q_\lambda\|_{C^\tau S^{-\alpha'}_{1,0;N}} \leq C_{\delta'}(1+|\lambda|)^{-\frac{\alpha-\alpha'}\alpha}.
  \end{eqnarray*}
  Since $a_\lambda(x,\xi,x)=0$, we can use Proposition~\ref{prop:RedOp} to conclude that
  \begin{equation*}
   a_\lambda(x,D_x,x)= \sum_{j=0}^\infty a_{\lambda,j}(x,D_x,x)
  \end{equation*}
  is well-defined as limit in $\mathcal{L}(\sC^s_0(\Rn))$ and satisfies
  \begin{equation*}
    \|a_\lambda(x,D_x,x)\|_{\mathcal{L}(\sC^s_0(\Rn))}\leq C \|a_\lambda\|_{C^\tau S^{\alpha-\alpha'}_{1,0;N}}\leq  C_{\delta'}(1+|\lambda|)^{-\frac{\alpha-\alpha'}\alpha}
  \end{equation*}
  for each $0<\alpha'< \alpha$ with $\alpha'<\tau-s$.
\end{proof}

Recall that an unbounded operator $A\colon \mathcal{D}(A)\subseteq X \to X$ generates an analytic semi-group on a Banach space $X$ if and only if $A$ is closed, $\mathcal{D}(A)$ is dense, and there are some $\delta>\frac{\pi}2$, $\omega\in\R$, and $M\geq 1$ such that $(\lambda-A)^{-1}$ exists for all $\lambda\in \omega+\Sigma_\delta$ and satisfies
\begin{equation}\label{eq:EstimGenerator}
  \|(\lambda-A)^{-1}\|_{\mathcal{L}(X)}\leq \frac{M}{|\lambda-\omega|}\qquad
\text{for all}\ \lambda\in \omega+\Sigma_\delta,
\end{equation}
cf. \cite{Pazy}. 

\begin{cor}\label{cor:Generator}
  Let $0<s<\tau$.
  Then $p(x,D_x)$ and $L$ generate an analytic semi-group on $\lh^s(\Rn)$ with domains $\mathcal{D}(L)=\mathcal{D}(p(x,D_x))= \lh^{s+\alpha}(\Rn)$.
Moreover, if $A=p(x,D_x)$ or $A=L$, then
  \begin{equation*}
    \|(\lambda-A)^{-1}\|_{\mathcal{L}(\sC^{s}_0(\Rn), \sC^{s+\alpha'}_0(\Rn))}
    \leq C_{\delta'}|\lambda|^{-\frac{\alpha-\alpha'}\alpha}\qquad
      \text{for all} \ \lambda\in \Sigma_{\delta`}, |\lambda|\geq R,
  \end{equation*}
  for all $0\leq \alpha'\leq \alpha$ with some sufficiently large $R>0$ and some $\delta'>\frac{\pi}2$.
\end{cor}
\begin{proof} 
  By a standard Neumann series argument 
  Lemma~\ref{lem:Parametrix} yields that 
$$
(\lambda-p(x,D_x))^{-1}\colon \sC^{s}_0(\Rn)\to \sC^{s+\alpha}_0(\Rn)
$$ 
exists for
  all $\lambda\in \Sigma_{\delta'}$ with $|\lambda|\geq R$ for some $R>0$ and
  satisfies
  \begin{equation*}
    \|(\lambda-p(x,D_x))^{-1}\|_{\mathcal{L}(\sC^s_0(\Rn))} \leq 2 \|q_\lambda(D_x,x)\|_{\mathcal{L}(\sC^s_0(\Rn))}\leq C |\lambda|^{-1}.
  \end{equation*}
  This implies (\ref{eq:EstimGenerator}) for a suitable choice of $\omega$.
  Hence $p(x,D_x)$ generates an analytic semi-group on $\lh^s(\Rn)$ with domain $\mathcal{D}(p(x,D_x))= \lh^{s+\alpha}(\Rn)$.

  Similarly,
  \begin{equation*}
    (\lambda - L) q_\lambda(D_x,x) = I - R_\lambda + L^2 q_\lambda (D_x,x),
  \end{equation*}
  where
  \begin{equation*}
    \|L^2 q_\lambda(D_x,x)\|_{\mathcal{L}(\sC^s_0(\Rn))}\leq 
    C \|q_\lambda(D_x,x)\|_{\mathcal{L}(\sC^s_0(\Rn),\sC^{s+\alpha''}_0(\Rn))}
    \leq C_{\delta,\delta',\alpha''}|\lambda|^{-\frac{\alpha-\alpha''}\alpha}
  \end{equation*}
  uniformly in $\lambda\in\Sigma_{\delta'}$, $|\lambda|\geq R$, with arbitrary $\alpha'<\alpha''<\alpha$. Thus the same arguments as before show that
  $L$ generates an analytic semi-group. 

  Finally, the uniform estimate of $(\lambda-A)^{-1}$ easily follows from Corollary~\ref{cor:MappingOfParametrix}.
\end{proof}

\noindent
\begin{proof*}{of Theorem~\ref{thm:main1}}
Because of Corollary~\ref{cor:Generator}, well-known results from semi-group theory imply the existence of a unique classical solution $u\in C^{1,\theta}([0,T];\lh^s(\Rn))\cap C^\theta([0,T];\mathcal{D}(L))$ of (\ref{eq:Parabol1})-(\ref{eq:Parabol2}), cf. \cite[Chapter 4, Theorem 3.5]{Pazy}. Finally, since $(\lambda-L)^{-1} \colon \lh^s(\Rn)\to \lh^{s+\alpha}(\Rn) $ is a bounded operator for $\lambda=R$, the graph norm on $\mathcal{D}(L)$, i.e., $\|u\|_{\sC^s}+\|Lu\|_{\sC^s}$, is equivalent to the norm of $\sC^{s+\alpha}(\Rn)$. That $u$ inherits the non-negativity from $f$ is easily established using the maximum principle. 
\end{proof*}

\section{The Martingale Problem}\label{sec:MP}

The standard reference for the martingale problem for diffusion operators is \cite{StVa79}. Since the paths of jump processes are not continuous by nature we have to set up the martingale problem for the path space $\sD([0,\infty);\R^n)$ of all c\`{a}dl\`{a}g paths. Good sources for this space are \cite{Bil99}, \cite{EtKu86}, \cite{JaSh03}, the first edition \cite{Bil68} is sufficient for many purposes. The standard reference for the martingale problem on $\sD([0,\infty);\R^n)$ is \cite{EtKu86}.

We denote by $\sD([0,\infty);\R^n)$ the set of all functions $\omega:[0,\infty)\to\R^n$ satisfying for all $t\geq 0$
\begin{align*}
\lim\limits_{s\to t+} \omega(s) = \omega(t) \,, \qquad \exists \, \omega(t-) = \lim\limits_{s\to t-} \omega(s) \,.
\end{align*}
A basic fact about $\sD([0,\infty);\R^n)$ is that any $\omega \in \sD([0,\infty);\R^n)$ has at most countably many points of discontinuity. As on the space of continuous functions the mapping $d_{uc}$ defined by
\[ d_{uc}(\omega_1, \omega_2) = \suml_{k\in\N} 2^{-k} \min\big\{1, \supl_{t\leq k} |\omega_1(t) - \omega_2(t)| \big\} \]
defines a metric. The space $\big(\sD([0,\infty);\R^n), d_{uc}\big)$ is a complete metric space but, different from the case of continuous functions, it is not separable. To see this, consider
\[ M:=\Big\{ \omega_s \in \sD([0,\infty);\R^n); \omega_s(t) = \mathbbm{1}_{[s,\infty)}(t), s \in[0,1)  \Big\} \,.\]
There cannot be a countable dense subset $A$ to the uncountable set $M$ since $d_{uc}(\omega_s, \omega_t)=\frac12$ as along as $s\ne t$. The set $A$ would need to be uncountable right away.

Nevertheless, there exists a metrizable topology on $\sD([0,\infty);\R^n)$ such that it becomes a complete, separable metric space. We summarize the main results on this space in the following theorem. Since the space $\sD([0,\infty);\R^n)$ is not too well known among analysts we include many details in this theorem. It is almost identical to Theorem VI.1.14 in \cite{JaSh03}.  

\begin{theorem} (1) There exists a metrizable topology on $\sD([0,\infty);\R^n)$, called the Skohorod topology for which the space is complete and separable. Denote the metric by $d$. Then $d(\omega_n, \omega) \to 0$ is equivalent to the existence of a sequence of strictly increasing functions $\lambda_n: [0,\infty) \to [0,\infty)$, satisfying $\lambda_n(0)=0, \lambda_n(t) \nearrow \infty$ for $t \to \infty$ and at the same time
\begin{align*}
\begin{cases}
&\supl_{s\geq 0} |\lambda_n(s) -s | \to 0 \quad \text{ as } n \to \infty \,, \\
&\Big(\supl_{s \leq k} |\omega_n (\lambda_n(s)) - \omega(s)| \to 0 \quad \text{ as } n \to \infty \Big) \quad \forall \, k \in \N \,. 
\end{cases}
\end{align*}
(2) A set $M \subset \sD([0,\infty);\R^n)$ is relatively compact for the Skohorod topology if and only if
\begin{align*}
\begin{cases}
&\supl_{\omega \in M} \supl_{s\leq k} |\omega(s)| < \infty \quad \forall k \in \N \,, \\
&\liml_{\rho \to 0+} \supl_{\omega \in M} \gamma_k(\omega,\rho) = 0 \quad \forall \, k \in \N \,. 
\end{cases}
\end{align*}
where $\gamma_k(\omega,t)$ is a generalized modulus of continuity, defined via
\begin{align*}
\gamma_k(\omega,\rho) = \inf\Big\{ \maxl_{i \leq L} \gamma(\omega;[t_{i-1},t_i)):0=t_0<\ldots<t_L=k, \infl_{i < L} (t_i - t_{i-1}) \geq \rho  \Big\} \,, 
\end{align*}
where $\gamma(\omega;I)$ is the usual modulus of continuity for $\omega$ on the interval $I \subset \R$. \\
(3) For given $t \geq 0$ let us denote by $\Pi_t$ the projection $\sD([0,\infty);\R^n)\to \R^n, \omega \mapsto \omega(t)=\Pi_t(\omega)$. With this notation the Borel $\sigma$-field $\sB\big(\sD([0,\infty);\R^n),d \big)$ equals $\sigma(\Pi_t; t \geq 0)$. \\
(4) The vector space $\big(\sD([0,\infty);\R^n),d\big)$ is not a topological vector space since addition of two elements is not continuous with respect to this topology.
\end{theorem}

A stochastic process $X$ with paths in $\sD([0,\infty);\R^n)$ can be interpreted as a random variable 
\[ X \colon (\Omega, \sF, \PP) \to \sD([0,\infty);\R^n) \]
with $X_t(\omega)=\omega(t)$ where $(\Omega, \sF, \PP)$ is an abstract probability space. Given a family $(X^{\alpha})_{\alpha \in \sA}$ of such processes we say that $(X^{\alpha})_{\alpha \in \sA}$ is relatively compact if the family $(\PP_{X^\alpha})_{\alpha \in \sA}$ of image measures $\PP_{X^\alpha}=\PP \circ (X^\alpha)^{-1}$ is relatively compact which, due to Prokhorov's theorem, amounts to saying that $(\PP_{X^\alpha})_{\alpha \in \sA}$ is tight. 

Usually, well-posedness of the martingale problem is much harder to be proved than mere solvability. A key feature that we use in order to show uniqueness of the solutions is formulated in the following lemma. It says that finite-dimensional distributions form a convergence determining class, see Theorem 3.7.8. in \cite{EtKu86}. 

\begin{lem}
Suppose that $(X^n)_{n \in \N}$ is a family of stochastic processes $X^n:(\Omega, \sF, \PP) \to \sD([0,\infty);\R^n)$ such as $X$ and there is a dense subset $J \subset [0,\infty)$ such that
\begin{align*}
\big( X^n(t_1), \ldots X^n(t_N) \big) \overset{d}{\Rightarrow} \big( X(t_1), \ldots X(t_N) \big) \\
\text{ or, equivalently } \qquad  \PP_{\big( X^n(t_1), \ldots X^n(t_N) \big)} \to \PP_{\big( X(t_1), \ldots X(t_N) \big)} \text{ weakly } 
\end{align*}
for all finite subsets $\{t_1, \ldots , t_N\} \subset J$. Then $X^n \overset{d}{\Rightarrow} X$ or, equivalently $\PP_{X^n} \to \PP_{X}$ weakly. 
\end{lem}

The situation turns out to be even better for solutions to the martingale problem. The following universal result says that even one-dimensional distributions determine the measure provided they agree for all initial distributions $\mu$, see Theorem 4.4.2 in \cite{EtKu86}. 

\begin{lem}\label{lem:unique-one-dim}
Consider the linear operator $(L,D(L))$ with $L$ defined as in $(\ref {eq:DefnL})$. Assume that for any initial distribution $\mu$ and any two corresponding solutions $\PP^\mu$, $\Q^\mu$ to the martingale problem \[ \PP^\mu_{\Pi_t} = \Q^\mu_{\Pi_t} \quad \forall \, t \geq 0 \,,\] then there exists at most one solution to the martingale problem for any initial distribution $\mu$.
\end{lem} 

The key to the proof is to show that regular conditional probabilities solve the martingale problem. Finally, we can prove Theorem \ref{thm:main2}.

\begin{proof*}{of Theorem~\ref{thm:main2}}
The existence of a solution $\PP^\mu$ for a given distribution $\mu$ on $\R^n$ has been established by several authors, see Theorem 2.2 in \cite{Str75}, Theorem IX.2.31 in \cite{JaSh03} and Theorem 3.2 in \cite{Hoh94}. Thus, we concentrate on the question of uniqueness. Uniqueness follows from solvability of the deterministic parabolic equation (\ref{eq:Parabol2}). This is the strategy worked out in \cite{StVa79} for the case of diffusions. We show how it works in our situation.

Assume that there are two solutions $\PP^\mu$, $\Q^\mu$ to the martingale problem for a given distribution $\mu$. A key step is to show that, for any $T>0$ the stochastic process $M=(M_t)_{t\in[0,T]}$ defined via
\begin{align}
M_t = v(t,\Pi_t) - \il_0^t \big( \frac{\partial}{\partial s} + L\big) v (s,\Pi_s) \, ds \label{eq:parab-martingale}
\end{align}
is a $\PP^\mu$-martingale and thus also a $\Q^\mu$-martingale for any function $v \in C^{1,\theta}([0,T];\sC^s _0(\R^n))$ $\cap$ $C^{\theta}([0,T];\sC^{s+\alpha}_0(\R^n))$ with $s,\theta \in (0,1)$.  This is proved exactly as in Theorem 4.2.1. (ii) of \cite{StVa79}. There is no need to have second spatial derivatives of $u$ since $L$ is an integro-differential operator of order $\alpha$. 

The main result follows once the following equality 
\begin{align}
\il_0^T \phi(s) \E_{\PP^\mu} \big(\psi(\Pi_s)\big) \, ds = \il_0^T \phi(s) \E_{\Q^\mu} \big(\psi(\Pi_s)\big) \, ds \label{eq:unique-one-dist-test}
\end{align}
is established for any $T>0$ and any choice of $\phi \in C^\infty_0((0,T))$, $\psi\in C^\infty_0(\R^n)$. Here, $\E_{\PP^\mu}$ and $\E_{\Q^\mu}$ denote the expectation with respect to $\PP^\mu$ and $\Q^\mu$ respectively. Equality (\ref{eq:unique-one-dist-test}) proves the equality of one-dimensional distributions, i.e. $\PP^\mu_{\Pi_t} = \Q^\mu_{\Pi_t}$ for all $t>0$, which in light of Lemma \ref{lem:unique-one-dim} proves the desired uniqueness result.  Equality (\ref{eq:unique-one-dist-test}) is proved as follows.

Setting $f(t,x)=\phi(t)\psi(x)$, Theorem \ref{thm:main1} proves that there is a function $v$ belonging to $C^{1,\theta}([0,T];\sC^s_0(\R^n))$ $\cap$ $C^{\theta}([0,T];\sC^{s+\alpha}_0(\R^n))$ and solving 
\begin{alignat*}{2}
    \partial_t v + Lv &= f &\qquad& \text{in}\ (0,T)\times \Rn,\\
    v(T,\cdot)&= 0 && \text{in}\ \Rn. 
\end{alignat*}
Thus
\begin{align*}
-\il_0^T \phi(s) & \E_{\PP^\mu} \big(\psi(\Pi_s)\big) \, ds = - \E_{\PP^\mu} \il_0^T f(s,\Pi_s) \, ds = \E_{\PP^\mu} \big(M_T\big) = \E_{\PP^\mu} \big(M_0\big) \\
&= \E_{\PP^\mu} \big(v(0,\Pi_0)\big) = \il_{\sD([0,\infty);\R^n)} v(0,\Pi_0(\omega)) \PP^\mu(d\omega) = \il_{\R^n} v(0,x) \mu(dx) \,. 
\end{align*}
Since the same line with the same right-hand side holds true when $\PP^\mu$ is replaced by $\Q^\mu$ equality (\ref{eq:unique-one-dist-test}) is established. The theorem is proved. 
\end{proof*}
\noindent


\def\polhk#1{\setbox0=\hbox{#1}{\ooalign{\hidewidth
  \lower1.5ex\hbox{`}\hidewidth\crcr\unhbox0}}}
  \def\polhk#1{\setbox0=\hbox{#1}{\ooalign{\hidewidth
  \lower1.5ex\hbox{`}\hidewidth\crcr\unhbox0}}}
  \def\polhk#1{\setbox0=\hbox{#1}{\ooalign{\hidewidth
  \lower1.5ex\hbox{`}\hidewidth\crcr\unhbox0}}} \def\cprime{$'$}


\medskip

\noindent{\bf Helmut Abels}\\
Max Planck Institute for Mathematics in the Sciences \\
Inselstra\ss e 22 \\
D-04103 Leipzig, Germany\\
{\it abels@mis.mpg.de}\\

\medskip

\noindent{\bf Moritz Kassmann}\\
Institut f\"{u}r Angewandte Mathematik\\
Universit\"{a}t Bonn \\
Beringstrasse 6\\
D-53115 Bonn, Germany\\
{\it kassmann@iam.uni-bonn.de}

\end{document}